\documentclass[review]{elsarticle}
\usepackage{tikz}
\usepackage{lineno,hyperref}
\pdfoutput=1
\modulolinenumbers[5]
\usepackage{float}
\RequirePackage{amsmath}[2016/11/05]
\usepackage[ruled,vlined]{algorithm2e}
\usepackage{natbib}
\journal{Computer Methods in Applied Mechanics and Engineering}
\date{}
\begin{document}

\begin{frontmatter}

\title{Residual-based adaptivity for two-phase flow simulation in porous media using Physics-informed Neural Networks}

\author[1,2,3]{John Hanna}
\author[1]{Jos\'e V. Aguado}
\author[2]{Sebastien~Comas-Cardona}
\author[3]{Ramzi Askri}
\author[1]{Domenico Borzacchiello}

\address[1]{Institut de Calcul Intensif (ICI) at Ecole Centrale de Nantes, 1 Rue de la Noë, 44300 Nantes, France}
\address[2]{Institut de Recherche en Génie Civil et Mécanique (GEM) at Ecole Centrale de Nantes, 1 Rue de la Noë, 44300 Nantes, France}
\address[3]{Institut de Recherche Technologique mutualisé (IRT) Jules Verne, Chem. du Chaffault, 44340 Bouguenais, France}

\begin{abstract}

This paper aims to provide a machine learning framework to simulate two-phase flow in porous media. The proposed algorithm is based on Physics-informed neural networks (PINN). A novel residual-based adaptive PINN is developed and compared with the residual-based adaptive refinement (RAR) method and with PINN with fixed collocation points. The proposed algorithm is expected to have great potential to be applied to different fields where adaptivity is needed. In this paper, we focus on the two-phase flow in porous media problem. We provide two numerical examples to show the effectiveness of the new algorithm. It is found that adaptivity is essential to capture moving flow fronts. We show how the results obtained through this approach are more accurate than using RAR method or PINN with fixed collocation points, while having a comparable computational cost.
\end{abstract}


\end{frontmatter}

%\linenumbers

\section{Introduction}

Multi-phase flow and transport phenomena in porous media appear in a variety of industrial applications, such as injection-based fiber-reinforced composites processes \cite{rtm_multi_phase_1, rtm_multi_phase_2}, exploitation of oil reservoirs \cite{oil_multi_phase_1, oil_multi_phase_2} and water resources management \cite{water_multiphase_1, water_multiphase_2}. That is why efficient stable modeling techniques are always required to understand the underlying physics of flow mechanisms.

Many numerical techniques have been developed to compute solutions to the partial differential equations governing these models. Grid-based classical methods such as finite element, finite volume, and finite difference are generally seen as the state of the art due to their numerical efficiency and stability \cite{grid_based_book}. To capture the flow front movement across the computational grid, these techniques are coupled with methods such as level-set \cite{level_sets}, volume of fluids \cite{vof} or phase-field method \cite{phase_field}. On the other hand, meshless and particle-based methods \cite{sph} have been proved to provide more natural ways of tracking the interface although they are generally regarded as having a lower accuracy and reduced stability \cite{sph_inaccurate_1, sph_inaccurate}.
Practical applications include flow monitoring and control \cite{rtm_control}, identification of flow parameters such as porosity and permeability \cite{permiability_identifi}, uncertainty quantification \cite{uncertainity_quant} and flow optimization \cite{rtm_optimization}.   

A successful approach for solving inverse problems relies on building surrogate parametric models for fast exploration of the state space of parameters \cite{surrogate_inverse_2, surrogate_inverse_1}.  These meta-models are specifically designed to be an optimal trade-off between accuracy and computational cost and can be effectively used for inverse problems when the minimization of the objective function requires multiple calls to the flow solver to perform parametric sweeps. Even so, when the dimension of the parametric space is high, the resulting complexity grows exponentially making some problems computationally intractable.

In recent years model order reduction techniques have been developed to tame this curse of dimensionality \cite{pod, mor}. These rely on the fundamental assumption that the solution of the parametric problem lies in a low-dimensional manifold of the original subspace where the solution approximation is sought. Learning the structure of this manifold is done through an offline training procedure minimizing the L2 distance from existing data, usually collected from multiple runs of a high-fidelity solver, and produces a low-rank basis that can be reused to represent the solution of new problems for unseen choices of the parameters. In practice, the choice of reduced-bases representation is equivalent to assuming a tensor format for the problem solution. Among the different choices the Canonical Polyadic (CP), Tucker, and Tensor Train (TT) are the most commonly used, as they provide a compact representation of the parametric solution as well as a reduced complexity of the model \cite{t_decompose_1, t_decompose_2, t_decompose_3}.   

The presence of a moving flow front in multi-phase flow introduces an additional difficulty to get effective model representations. Tensor formats are generally regarded as unfit to represent solutions exhibiting a moving discontinuity because, due to the dual-scale nature of the problem, there is a need for a large number of basis vectors to obtain a good approximation of the solution. For instance, applying a standard CP space-time decomposition to a simple 1D moving Heaviside function cannot provide an accurate approximation unless a high number of modes is used. In practice, the rank needed to obtain reasonably accurate results is not offering any computational advantage compared to full-order representations. This issue is well documented in the community of model order reduction and it affects not only hyperbolic equations that are likely to give rise to shocks or discontinuities but also to other problems in which the physics involved produces localized effects in the solutions. The issue was tackled by \cite{mor_isssue_1, mor_issue_2, mor_issue_3} who proposed ways to fix the problem, however, it is still a pressing difficulty that requires attention. 

In recent years, with the ever-increasing number of machine learning applications that have emerged, neural networks (NN) have established themselves as a class of architectures offering superior modeling capabilities in problems such as regression and classification. When compared to the previously mentioned tensor formats NNs might provide a more flexible tool to represent functions with moving discontinuities resulting from the solution of transport equations.

It is the main aim of this paper to assess this opportunity and provide the basis for a proper framework to solve multi-phase flow in porous media through the use of NNs.

The proposed numerical technique is based on Physics-informed neural networks (PINN), a PDE solver based on the use of neural networks approximation as the search space for the PDE solution. A solution to a PDE can be obtained by minimizing the PDE and boundary conditions residuals over a finite set of points called collocation and training points, respectively. PINN was first introduced in \cite{pinn_main} and has been used to solve many forward and inverse problems \cite{xpinn, conservative_pinn, high_speed_pinn, fractional_pinn, ultrasound_pinn, additive_pinn, solid_pinn, composite_pinn, surrogate_pinn, subsurface_pinn, fracture_pinn}. We show in this paper that the choice of collocations points is crucial to obtain accurate results.

The main contributions of this paper can be summarized as:

\begin{itemize}
    \item Providing a basic framework to simulate two-phase flow in porous media using PINN.
    \item Developing a residual-based adaptive PINN for accurate flow front predictions.
    \item Extending the adaptivity algorithm for the training points to better capture the initial/boundary conditions and to have different collocation points for different PDEs in a coupled system.
\end{itemize}

The rest of the paper is organized as follows. In section \ref{section:pinn}, the PINN method to solve two-phase flow in porous media is described. The novel residual-based adaptivity algorithm is detailed in section \ref{section:adaptivity}. Numerical experiments were carried out and presented in section \ref{section:results}, showing the advantages of using the adaptivity technique introduced in this work compared to using classical PINN with fixed collocation points or the RAR technique. Finally, a summary and a conclusion are given in section \ref{section:conclusion}.

\section{Model problem and PINN formulation}
\label{section:pinn}

In this section, we first introduce the model problem for two-phase flow in porous media. The basics for PINN are then recalled. Finally, PINN is applied to solve the model problem.

 \subsection{Model problem for two-phase flow in porous media}

Flow in porous media can be described by Darcy's law, which reads as

\begin{equation}
    \mathbf{v} = -\dfrac{1}{\mu}\mathbf{K}\cdot\nabla p 
    \label{eq:darcy_full}
\end{equation}

where $\mathbf{v}$ is the volume average Darcy's velocity, $\mathbf{K}$ the permeability tensor, $\mu$ the viscosity, and $\nabla p$ the pressure gradient. Both fluids are assumed to be incompressible, therefore, the mass conservation equation reduces to

\begin{equation}
   \nabla \cdot \mathbf{v} = 0
    \label{eq:momentum_full}
\end{equation}

To differentiate between the two fluid phases, a fraction function $c$ is introduced which takes a value $1$ for one fluid and $0$ for the other one. The viscosity $\mu$ is redefined as

\begin{equation}
    \mu = c\mu_2 + (1-c)\mu_1
    \label{eq:mu}
\end{equation}

where $\mu_2$ and $\mu_1$ are the viscosities of the two fluids. $c$ evolves with time according to the following advection equation

\begin{equation}
    c_t + \mathbf{v}\cdot\nabla c = 0
    \label{eq:c_full}
\end{equation}

where $c_t$ is the time derivative of the fraction function $c$.

The problem formulation is completed by assigning boundary conditions for  $\mathbf{v}$, $p$ and $c$, as well as initial conditions for $c$:

\begin{equation}\label{eq:c_IC}
    c(\mathbf{x},t=0)=c_0(\mathbf{x})\, .
\end{equation}

Pressure Dirichlet boundary conditions are prescribed on inlet and outlet boundaries:

\begin{equation}\label{eq:p_BC_in}
    p(\mathbf{x}_{inlet},t) = p_{in} \quad \text{(Inlet pressure)}
\end{equation}

and 

\begin{equation}\label{eq:p_BC_out}
    p(\mathbf{x}_{outlet},t) = p_{out} \quad  \text{(Outlet pressure)}
\end{equation}

Impermeable boundaries are characterized by zero normal velocity:

\begin{equation}\label{eq:v_BC}
    \mathbf{v\cdot n}  = 0  \quad \text{(Impermeable wall)}
\end{equation}

Inlet flow also requires the assignment of boundary conditions for $c$:

\begin{equation}\label{eq:c_BC}
    c(\mathbf{x}_{inlet},t) = 1 \quad \text{(Inlet)}
\end{equation}

\subsection{PINN}

The first basic idea in PINN is the choice of the search space as a fully connected neural network approximation. The approximation has the form of

\begin{equation}
    \mathbf{z}^i = \sigma^i(\mathbf{W}^i \mathbf{z}^{i-1} + \mathbf{b}^i),\ \ i=1,...,L
\end{equation}

where $\mathbf{z^0}$ and $\mathbf{z^L}$ are the input (temporal or space dimensions) and output (an approximation to the solution) of the network, respectively. While, $\mathbf{W^i}$ and $\mathbf{b^i}$ are the parameters of each layer, known as the weights and biases, respectively. $\sigma^i$ is a nonlinear function use to add nonlinearity in the representation and is called the activation function. $L$ is the number of layers in the network.

Secondly, the residual of the PDE, to be solved, is obtained by differentiating the neural network output with respect to the inputs using automatic differentiation. A solution is reached by finding the weights and biases that minimize a loss function, composed of the PDEs and boundary conditions residual over a set of finite points namely, collocation and training points, respectively.

\subsection{PINN structure and loss function definition for the model problem}

\begin{figure}[h]
    \centering
\tikzset{every picture/.style={line width=0.75pt}} %set default line width to 0.75pt        

\begin{tikzpicture}[x=0.75pt,y=0.75pt,yscale=-1,xscale=1]

%Shape: Circle [id:dp9677571531571376] 
\draw   (0,85) .. controls (0,71.19) and (11.19,60) .. (25,60) .. controls (38.81,60) and (50,71.19) .. (50,85) .. controls (50,98.81) and (38.81,110) .. (25,110) .. controls (11.19,110) and (0,98.81) .. (0,85) -- cycle ;
%Shape: Circle [id:dp16473337849285352] 
\draw   (3,168) .. controls (3,154.19) and (14.19,143) .. (28,143) .. controls (41.81,143) and (53,154.19) .. (53,168) .. controls (53,181.81) and (41.81,193) .. (28,193) .. controls (14.19,193) and (3,181.81) .. (3,168) -- cycle ;
%Shape: Circle [id:dp7181373137331066] 
\draw   (300.8,48) .. controls (300.8,34.19) and (311.99,23) .. (325.8,23) .. controls (339.61,23) and (350.8,34.19) .. (350.8,48) .. controls (350.8,61.81) and (339.61,73) .. (325.8,73) .. controls (311.99,73) and (300.8,61.81) .. (300.8,48) -- cycle ;
%Shape: Circle [id:dp34971532434167996] 
\draw   (303,129) .. controls (303,115.19) and (314.19,104) .. (328,104) .. controls (341.81,104) and (353,115.19) .. (353,129) .. controls (353,142.81) and (341.81,154) .. (328,154) .. controls (314.19,154) and (303,142.81) .. (303,129) -- cycle ;
%Shape: Circle [id:dp22064826925129766] 
\draw   (302.8,209) .. controls (302.8,195.19) and (313.99,184) .. (327.8,184) .. controls (341.61,184) and (352.8,195.19) .. (352.8,209) .. controls (352.8,222.81) and (341.61,234) .. (327.8,234) .. controls (313.99,234) and (302.8,222.81) .. (302.8,209) -- cycle ;
%Shape: Rectangle [id:dp06478446043740993] 
\draw   (149,28) -- (246.8,28) -- (246.8,68) -- (149,68) -- cycle ;
%Shape: Rectangle [id:dp08485878916036715] 
\draw   (389.8,23) -- (560.8,23) -- (560.8,63) -- (389.8,63) -- cycle ;
%Straight Lines [id:da41186878537511196] 
\draw    (50,85) -- (144.94,47.73) ;
\draw [shift={(146.8,47)}, rotate = 158.57] [color={rgb, 255:red, 0; green, 0; blue, 0 }  ][line width=0.75]    (10.93,-3.29) .. controls (6.95,-1.4) and (3.31,-0.3) .. (0,0) .. controls (3.31,0.3) and (6.95,1.4) .. (10.93,3.29)   ;
%Straight Lines [id:da7060621644195328] 
\draw    (53,168) -- (145.57,48.58) ;
\draw [shift={(146.8,47)}, rotate = 127.78] [color={rgb, 255:red, 0; green, 0; blue, 0 }  ][line width=0.75]    (10.93,-3.29) .. controls (6.95,-1.4) and (3.31,-0.3) .. (0,0) .. controls (3.31,0.3) and (6.95,1.4) .. (10.93,3.29)   ;
%Straight Lines [id:da2170515368242385] 
\draw    (50,85) -- (147.97,128.39) ;
\draw [shift={(149.8,129.2)}, rotate = 203.89] [color={rgb, 255:red, 0; green, 0; blue, 0 }  ][line width=0.75]    (10.93,-3.29) .. controls (6.95,-1.4) and (3.31,-0.3) .. (0,0) .. controls (3.31,0.3) and (6.95,1.4) .. (10.93,3.29)   ;
%Straight Lines [id:da725832785396104] 
\draw    (50,85) -- (149.55,209.64) ;
\draw [shift={(150.8,211.2)}, rotate = 231.38] [color={rgb, 255:red, 0; green, 0; blue, 0 }  ][line width=0.75]    (10.93,-3.29) .. controls (6.95,-1.4) and (3.31,-0.3) .. (0,0) .. controls (3.31,0.3) and (6.95,1.4) .. (10.93,3.29)   ;
%Straight Lines [id:da048711002907182] 
\draw    (53,168) -- (147.94,129.94) ;
\draw [shift={(149.8,129.2)}, rotate = 158.16] [color={rgb, 255:red, 0; green, 0; blue, 0 }  ][line width=0.75]    (10.93,-3.29) .. controls (6.95,-1.4) and (3.31,-0.3) .. (0,0) .. controls (3.31,0.3) and (6.95,1.4) .. (10.93,3.29)   ;
%Straight Lines [id:da5496911076200688] 
\draw    (53,168) -- (148.97,210.39) ;
\draw [shift={(150.8,211.2)}, rotate = 203.83] [color={rgb, 255:red, 0; green, 0; blue, 0 }  ][line width=0.75]    (10.93,-3.29) .. controls (6.95,-1.4) and (3.31,-0.3) .. (0,0) .. controls (3.31,0.3) and (6.95,1.4) .. (10.93,3.29)   ;
%Straight Lines [id:da8003751394098775] 
\draw    (246.8,48) -- (299.53,47.81) ;
\draw [shift={(301.53,47.8)}, rotate = 179.79] [color={rgb, 255:red, 0; green, 0; blue, 0 }  ][line width=0.75]    (10.93,-3.29) .. controls (6.95,-1.4) and (3.31,-0.3) .. (0,0) .. controls (3.31,0.3) and (6.95,1.4) .. (10.93,3.29)   ;
%Straight Lines [id:da3415989983746275] 
\draw    (246.8,129) -- (300.53,128.81) ;
\draw [shift={(302.53,128.8)}, rotate = 179.79] [color={rgb, 255:red, 0; green, 0; blue, 0 }  ][line width=0.75]    (10.93,-3.29) .. controls (6.95,-1.4) and (3.31,-0.3) .. (0,0) .. controls (3.31,0.3) and (6.95,1.4) .. (10.93,3.29)   ;
%Straight Lines [id:da04868931690134737] 
\draw    (248.8,210) -- (298.53,210.77) ;
\draw [shift={(300.53,210.8)}, rotate = 180.89] [color={rgb, 255:red, 0; green, 0; blue, 0 }  ][line width=0.75]    (10.93,-3.29) .. controls (6.95,-1.4) and (3.31,-0.3) .. (0,0) .. controls (3.31,0.3) and (6.95,1.4) .. (10.93,3.29)   ;
%Shape: Rectangle [id:dp22556842915771758] 
\draw   (151,189) -- (248.8,189) -- (248.8,229) -- (151,229) -- cycle ;
%Shape: Rectangle [id:dp7037424331232982] 
\draw   (148,110) -- (245.8,110) -- (245.8,150) -- (148,150) -- cycle ;
%Shape: Rectangle [id:dp3975649620949018] 
\draw   (392,106) -- (559.8,106) -- (559.8,146) -- (392,146) -- cycle ;
%Shape: Rectangle [id:dp39622132877259264] 
\draw   (393,190) -- (559.8,190) -- (559.8,230) -- (393,230) -- cycle ;
%Straight Lines [id:da6924279131629054] 
\draw    (350.8,48) -- (387.83,41.35) ;
\draw [shift={(389.8,41)}, rotate = 169.82] [color={rgb, 255:red, 0; green, 0; blue, 0 }  ][line width=0.75]    (10.93,-3.29) .. controls (6.95,-1.4) and (3.31,-0.3) .. (0,0) .. controls (3.31,0.3) and (6.95,1.4) .. (10.93,3.29)   ;
%Straight Lines [id:da901530805461241] 
\draw    (352.8,209) -- (389.37,42.95) ;
\draw [shift={(389.8,41)}, rotate = 102.42] [color={rgb, 255:red, 0; green, 0; blue, 0 }  ][line width=0.75]    (10.93,-3.29) .. controls (6.95,-1.4) and (3.31,-0.3) .. (0,0) .. controls (3.31,0.3) and (6.95,1.4) .. (10.93,3.29)   ;
%Straight Lines [id:da8688473487547783] 
\draw    (352.8,209) -- (389.96,128.81) ;
\draw [shift={(390.8,127)}, rotate = 114.86] [color={rgb, 255:red, 0; green, 0; blue, 0 }  ][line width=0.75]    (10.93,-3.29) .. controls (6.95,-1.4) and (3.31,-0.3) .. (0,0) .. controls (3.31,0.3) and (6.95,1.4) .. (10.93,3.29)   ;
%Straight Lines [id:da5209265065511064] 
\draw    (350.8,48) -- (389.9,125.22) ;
\draw [shift={(390.8,127)}, rotate = 243.15] [color={rgb, 255:red, 0; green, 0; blue, 0 }  ][line width=0.75]    (10.93,-3.29) .. controls (6.95,-1.4) and (3.31,-0.3) .. (0,0) .. controls (3.31,0.3) and (6.95,1.4) .. (10.93,3.29)   ;
%Straight Lines [id:da9522841641779562] 
\draw    (350.8,48) -- (392.31,214.06) ;
\draw [shift={(392.8,216)}, rotate = 255.96] [color={rgb, 255:red, 0; green, 0; blue, 0 }  ][line width=0.75]    (10.93,-3.29) .. controls (6.95,-1.4) and (3.31,-0.3) .. (0,0) .. controls (3.31,0.3) and (6.95,1.4) .. (10.93,3.29)   ;
%Straight Lines [id:da01514571468971293] 
\draw    (353,129) -- (388.8,127.11) ;
\draw [shift={(390.8,127)}, rotate = 176.97] [color={rgb, 255:red, 0; green, 0; blue, 0 }  ][line width=0.75]    (10.93,-3.29) .. controls (6.95,-1.4) and (3.31,-0.3) .. (0,0) .. controls (3.31,0.3) and (6.95,1.4) .. (10.93,3.29)   ;

% Text Node
\draw (14,74) node [anchor=north west][inner sep=0.75pt]  [font=\huge] [align=left] {$\displaystyle \mathbf{x}$};
% Text Node
\draw (22,162) node [anchor=north west][inner sep=0.75pt]  [font=\LARGE] [align=left] {$\displaystyle t$};
% Text Node
\draw (316,39) node [anchor=north west][inner sep=0.75pt]  [font=\LARGE] [align=left] {$\displaystyle \mathbf{v}$};
% Text Node
\draw (318,119) node [anchor=north west][inner sep=0.75pt]  [font=\LARGE] [align=left] {$\displaystyle p$};
% Text Node
\draw (319,201) node [anchor=north west][inner sep=0.75pt]  [font=\LARGE] [align=left] {$\displaystyle c$};
% Text Node
\draw (171,38) node [anchor=north west][inner sep=0.75pt]  [font=\large] [align=left] {1st NN};
% Text Node
\draw (169,118) node [anchor=north west][inner sep=0.75pt]  [font=\large] [align=left] {2nd NN};
% Text Node
\draw (172,198) node [anchor=north west][inner sep=0.75pt]  [font=\large] [align=left] {3rd NN};
% Text Node
\draw (416,31) node [anchor=north west][inner sep=0.75pt]  [font=\large]  {$f_{1} :=c_{t} +\mathbf{v} \cdot \nabla c\ $};
% Text Node
\draw (399.2,112) node [anchor=north west][inner sep=0.75pt]  [font=\large]  {$\mathbf{f_{2}} :=\mathbf{v} +( 1/\mu ) \ \mathbf{K} \cdot \nabla p$};
% Text Node
\draw (437.2,197) node [anchor=north west][inner sep=0.75pt]  [font=\large]  {$f_{3} :=\nabla \cdot \mathbf{v}$};
\end{tikzpicture}

    \caption{PINN structure for a general two-phase flow in porous media problem.}\label{fig:pinn_1d_fill}
    
\end{figure}
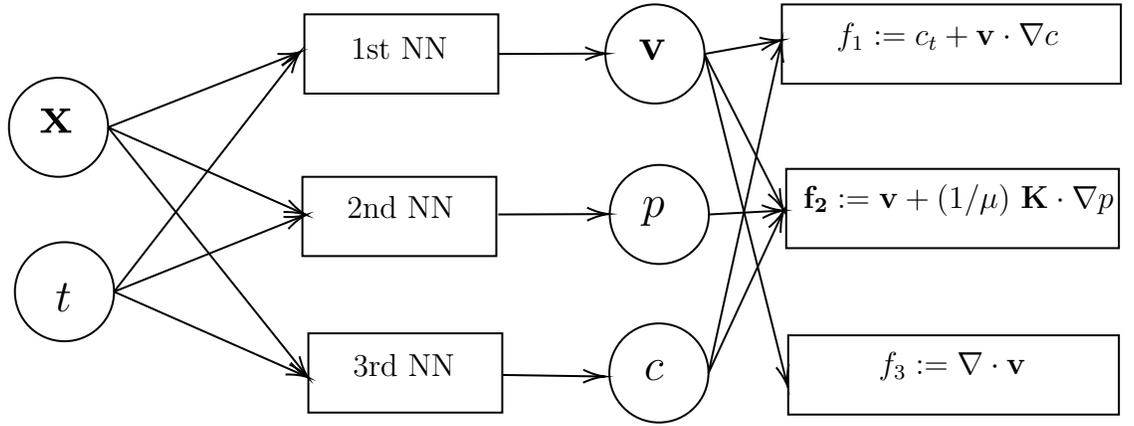

To solve the problem using PINN, three distinct neural networks are used: one for the velocity, another for pressure, and another for the fraction function. Each of these networks has space $\mathbf{x}$ and time $t$ inputs. The outputs of these neural networks ($\mathbf{v}$, $p$ and $c$) are differentiated with respect to the inputs, using automatic differentiation, forming the residuals of the three differential equations (\ref{eq:c_full}), (\ref{eq:darcy_full}) and (\ref{eq:momentum_full}). It should be noted that the network of $\mathbf{v}$ has 1, 2, or 3 outputs,  according to the problem dimension, corresponding to the velocity components. Another possible structure is to have different networks for each velocity component. The PINN structure for a general two-phase flow in porous media problem is summarized in figure \ref{fig:pinn_1d_fill}.

The loss function can, then, be defined as follows

\begin{equation}
    Loss = \lambda_v\ loss_v + \lambda_c\ loss_c +\lambda_p\ loss_p +\lambda_{1}\ loss_{f_1} +\lambda_{2}\ loss_{f_2} +\lambda_{3}\ loss_{f_3}
\end{equation}

where

\begin{gather}
    loss_v = \frac{1}{N_v}\sum_{i=1}^{N_v} r_v^2(\mathbf{x_v^i}, t_v^i) =\frac{1}{N_v}\sum_{i=1}^{N_v} ||\mathbf{v}(\mathbf{x_v^i}, t_v^i)\cdot \mathbf{n}(\mathbf{x_v^i}, t_v^i)||^2\\[2ex]
    loss_c = \frac{1}{N_c}\sum_{i=1}^{N_c} r_c^2(\mathbf{x_c^i}, t_c^i) =\frac{1}{N_c}\sum_{i=1}^{N_c} ||c(\mathbf{x_c^i}, t_c^i) - c_b^i||^2\\[2ex]
    loss_p = \frac{1}{N_p}\sum_{i=1}^{N_p} r_p^2(\mathbf{x_p^i}, t_p^i) = \frac{1}{N_p}\sum_{i=1}^{N_p} ||p(\mathbf{x_p^i}, t_p^i) - p_b^i||^2\\[2ex]
    loss_{f1} = \frac{1}{N_{f1}}\sum_{i=1}^{N_{f1}} ||f_1(\mathbf{x_{f1}^i}, t_{f1}^i)||^2= \frac{1}{N_{f1}}\sum_{i=1}^{N_{f1}} ||c_t + \mathbf{v}\cdot\nabla c||_{(\mathbf{x_{f1}^i}, t_{f1}^i)}^2\\[2ex]
    loss_{f2} = \frac{1}{N_{f2}}\sum_{i=1}^{N_{f2}} ||f_2(\mathbf{x_{f2}^i}, t_{f2}^i)||^2= \frac{1}{N_{f2}}\sum_{i=1}^{N_{f2}} ||\mathbf{v} + \frac{1}{\mu}\mathbf{K}\cdot\nabla p||_{(\mathbf{x_{f2}^i}, t_{f2}^i)}^2\\[2ex]
    loss_{f3} = \frac{1}{N_{f3}}\sum_{i=1}^{N_{f3}} ||f_3(\mathbf{x_{f3}^i}, t_{f3}^i)||^2= \frac{1}{N_{f3}}\sum_{i=1}^{N_{f3}} ||\nabla \cdot \mathbf{v}||_{(\mathbf{x_{f3}^i}, t_{f3}^i)}^2
\end{gather}

where $\{\mathbf{x_v^i},t_v^i,v_b^i\}_{i=1}^{N_v}$, $\{\mathbf{x_c^i},t_c^i,c_b^i\}_{i=1}^{N_c}$ and $\{\mathbf{\mathbf{x_p^i}},t_p^i,p_b^i\}_{i=1}^{N_p}$ are the points where the initial/boundary conditions are defined for $\mathbf{v}$, $c$ and $p$, respectively. While $\{\mathbf{x_{f1}^i},t_{f1}^i\}_{i=1}^{N_{f1}}$, $\{\mathbf{x_{f2}^i},t_{f2}^i\}_{i=1}^{N_{f2}}$ and $\{\mathbf{x_{f3}^i},t_{f3}^i\}_{i=1}^{N_{f3}}$ are the collocation points in space and time for the three residuals, respectively, where the physics is enforced and $\lambda_i$ are the weights of each term in the loss function. The $\lambda_i$ weights are important to make each contribution to the loss function has comparable magnitude, thus, aiding the optimization process. Finally, a solution for the fields, $p$, $\mathbf{v}$ and $c$, is obtained by minimizing the loss function with respect to the neural networks' parameters.

\section{Collocation points adaptivity}
\label{section:adaptivity}

Mesh Refinement is a basic idea in classical numerical techniques such as finite element (FEM) and finite volume (FVM) methods \citep{mesh_refine}. There are three basic techniques to mesh adaptation: h-adaptivity \cite{h_adaptivity}, r-adaptivity \cite{r_adaptivity} and p-adaptivity \cite{p_adaptivity}. H-adaptivity adds more nodes, thus increasing the degrees of freedom and the mesh connectivity. While r-adaptivity keeps the same number of nodes and degrees of freedom however the nodes are relocated while keeping the same connectivity. Finally, p-adaptivity increases the polynomial degrees of elements while keeping the mesh fixed. Other adaptivity methods exist that combines some of the basic techniques together such as: hp-adaptivity \cite{hp_adaptivity} and rh-adaptivity \cite{p_adaptivity}

There are mainly three drivers for mesh adaptation: error \citep{error_adapt_1}, PDE residual \citep{res_adapt} and solution features \citep{sol_adapt_2}.

Error-based adaptation adds more degrees of freedom where the solution error is high. This technique ignores the fact that the error is transported in the domain. Therefore, adapting where the error is high might ignore the region of the error source itself, where adaption is more useful \citep{error_adapt_bad}.

Residual-based adaptation refines the mesh where the discretized PDE residual is high. The residual can be seen as the source of error in the solution \citep{error_adapt_2}. Thus, refining where the residual is high is seen as a way of refining where the error source is. Therefore, this technique usually performs better than error-based adaptation. 

Solution-based adaptation utilizes the solution features such as gradients or discontinuities for adaptation. The philosophy behind this technique is that by using more points in these locations, these features can be resolved, thus, leading to improving the overall accuracy of the solution. However, if multiple features are present in the same problem, the adaptation results in over-refining some features while others are ignored. An example of this adaptivity failure can be found in \cite{solution_based_fail}.

The adaptation process is usually computationally expensive since certain requirements have to be satisfied and the mesh connectivity needs to be updated. Moreover, parallelization becomes complex for unstructured grids. However, in the case of PINN, changing the collocation points locations or adding more points are cheap processes. The main reason is that PINN is a meshless method, thus, there are no element volumes to take care of or mesh connectivity to update. Moreover, the approximation of the derivatives is independent of the collocation point position. Therefore, there is no discretization error resulting from the distribution of the collocation points. The only thing to do is identify the location where more points are needed.

\begin{figure}[H]
    \centering
    \includegraphics[width=1\textwidth]{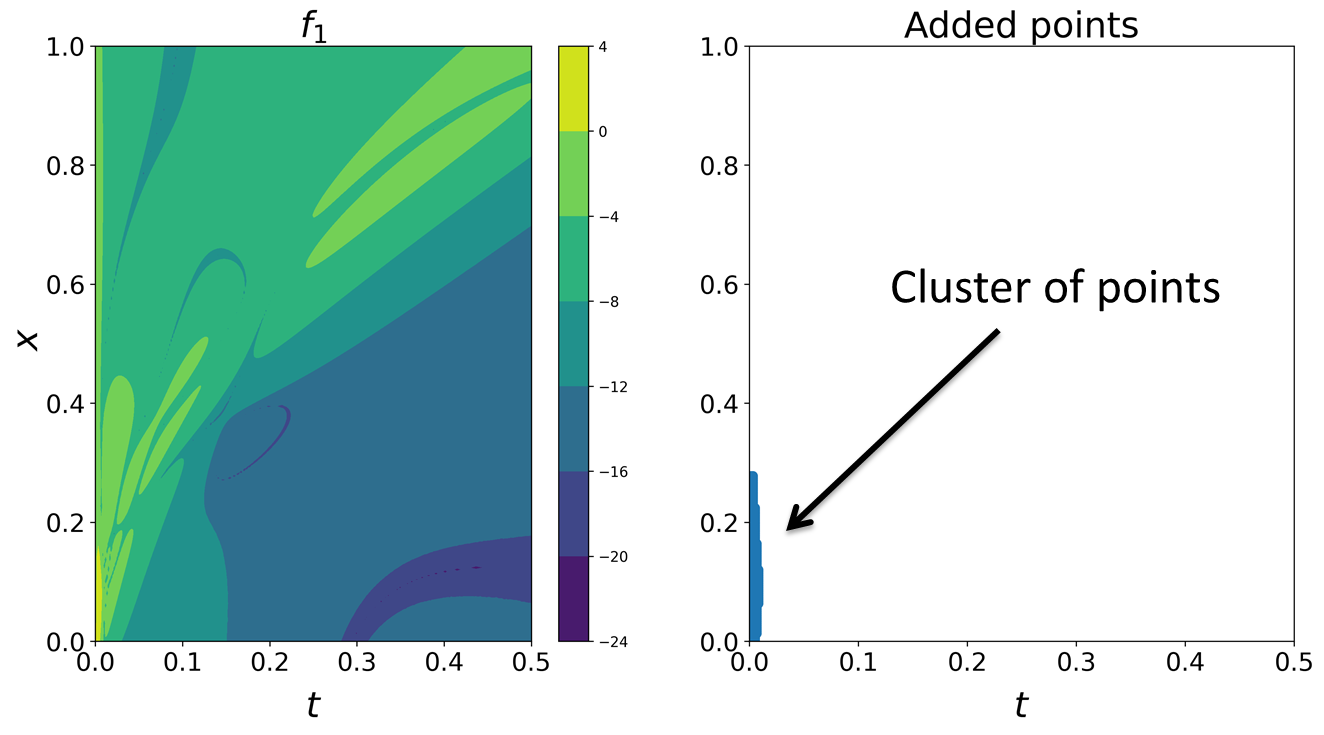}
  \caption{Left: log of the absolute residual field $f_1$, right: the chosen points according to RAR.}\label{fig:new_rar}
\end{figure}

In this study, we develop a residual-based algorithm by enriching the locations where the residual is high with more collocation points. We build the algorithm based on the work of \cite{res_adapt_pinn} in which the authors developed the residual-based adaptive refinement method (RAR). In their work, the authors used a dense set of randomly drawn points in the space-time domain, where residuals are evaluated. Points corresponding to the largest residual values are then added to the training set of collocation points. The progressive refinement of the training set allows for residual control. However when residual is showing high values in very narrow regions, this sampling strategy tends to produce excessively clustered points ignoring other solution features if existing and leading to unnecessary over-refinement. This behaviour appears to be related to solutions exhibiting moving sharp fronts or discontinuities, as in the model considered in this work. An example of that is seen in figure~\ref{fig:new_rar} where the added points are focused in a very small region.

\begin{figure}[H]
    \centering
    \includegraphics[width=1\textwidth]{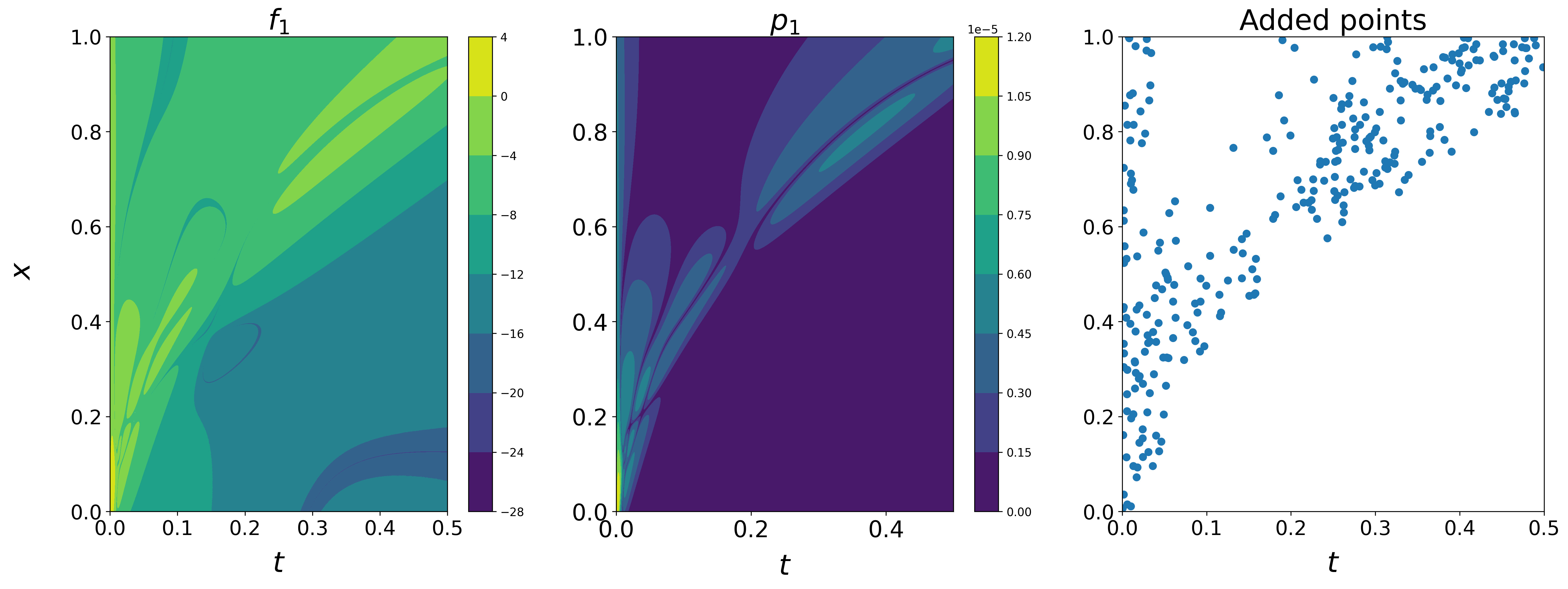}
  \caption{Left: log of the absolute residual field $f_1$, middle: the probability density function built from the residual, right: the chosen points drawn from the density function.}
    \label{fig:adaptivity}
\end{figure}

To avoid this point clusters that might lead to over-fitting of the model, we designed a probability density function based on the residual field to control the spread of the added points. The density function is, then, used to draw points from the dense set and these points are added to the training set. By doing this, collocation points will be more evenly spread in the domain; more points will be added where the probability is high (high residual) and fewer points where the probability is lower (low residual). Figure~\ref{fig:adaptivity} explains the procedure used for the adaptivity technique developed in this paper.

The algorithm is extended for coupled differential equations so that each PDE residual will have different collocation points. Moreover, data points are also enriched using a similar strategy to capture the initial/boundary conditions better. The algorithm is as follows:

\begin{algorithm}[H]
\SetAlgoLined
 Inputs: Number of adaptivity steps $M$, number of iterations $n$, tolerances $\epsilon_i$, $\epsilon_v$, $\epsilon_c$ and $\epsilon_p$\;
 \While{$m < M\ \&$ ($\mu_i > \epsilon_i\ ||\ \mu_v > \epsilon_v\ ||\ \mu_c > \epsilon_c\ ||\ \mu_p > \epsilon_p$) }{
  - Calculate $f_i$ (3 PDE residuals), $r_v$, $r_c$ and $r_p$ using dense sets\;
  - Build the probability functions $p_i$, $p_v$, $p_c$ and $p_p$ using dense sets Eq:~$\ref{eq:probability}$\;
  - Draw points from the dense sets using the probabilities and add it to the training sets\;
  - Use minimization algorithm for $n$ fixed iterations\;
  - Calculate the residuals' mean using the dense sets\;
   \begin{equation*}
    \begin{gathered}
        \mu_i = \frac{1}{N_i}\sum|f_i|\\
        \mu_v = \frac{1}{N_v}\sum|\mathbf{v}(\mathbf{x_v^i}, t_v^i)\cdot \mathbf{n}(\mathbf{x_v^i}, t_v^i)|\\
        \mu_c = \frac{1}{N_c}\sum|c(x_c^i,t_c^i)-c_b^i|\\
        \mu_p = \frac{1}{N_p}\sum|p(x_p^i,t_p^i)-p_b^i|
   \end{gathered}      
   \end{equation*}
 }
 \caption{Residual-based adaptivity}
\end{algorithm}

The probability functions used have the form of 

\begin{equation}
    p(\mathbf{X}) = \dfrac{max(log|r(\mathbf{X})/\epsilon|, 0)}{\int_{\Omega} max(log|r(\mathbf{X})/\epsilon|, 0) \,d\mathbf{X}}
    \label{eq:probability}
\end{equation}

where $\mathbf{X}$ is the random vector $[x, t]^T$, $r$ the considered residual, $\Omega$ the spatio-temporal domain and $\epsilon$  a small tolerance to filter small residual values. In practice, the choice of the value of $\epsilon$ is chosen to control the spread of the point distributions. The function is designed in a way to ensure that its integral over the space-time domain is 1, hence, the presence of the term in the denominator which is calculated using Monte Carlo integration over the dense set of points.

\section{Numerical experiments}\label{section:results}
\subsection{One dimensional injection}\label{subsec:1DFilling}
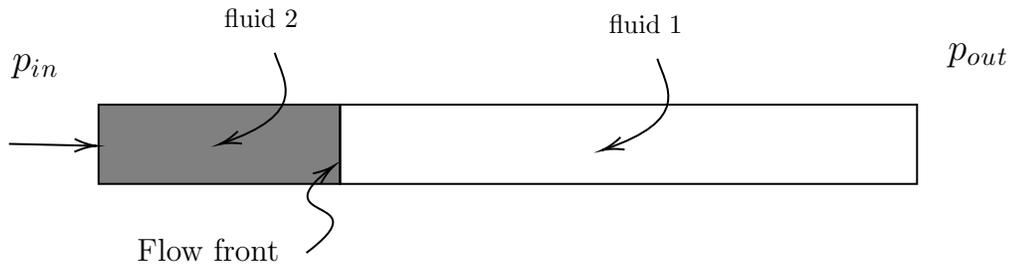
\begin{figure}[h]
    \centering

\tikzset{every picture/.style={line width=0.75pt}} %set default line width to 0.75pt        

\begin{tikzpicture}[x=0.75pt,y=0.75pt,yscale=-1,xscale=1]

%Shape: Rectangle [id:dp9935858859030073] 
\draw  [fill={rgb, 255:red, 128; green, 128; blue, 128 }  ,fill opacity=1 ] (123,93) -- (244.8,93) -- (244.8,133) -- (123,133) -- cycle ;
%Shape: Rectangle [id:dp8329373450379818] 
\draw   (244.8,93) -- (535.8,93) -- (535.8,133) -- (244.8,133) -- cycle ;
%Straight Lines [id:da3245034610645785] 
\draw    (77.8,112.8) -- (119.8,113.75) ;
\draw [shift={(121.8,113.8)}, rotate = 181.3] [color={rgb, 255:red, 0; green, 0; blue, 0 }  ][line width=0.75]    (10.93,-3.29) .. controls (6.95,-1.4) and (3.31,-0.3) .. (0,0) .. controls (3.31,0.3) and (6.95,1.4) .. (10.93,3.29)   ;
%Curve Lines [id:da2785954902937213] 
\draw    (227.8,167.8) .. controls (267.4,138.1) and (203.11,153.48) .. (240.63,124.69) ;
\draw [shift={(241.8,123.8)}, rotate = 503.13] [color={rgb, 255:red, 0; green, 0; blue, 0 }  ][line width=0.75]    (10.93,-3.29) .. controls (6.95,-1.4) and (3.31,-0.3) .. (0,0) .. controls (3.31,0.3) and (6.95,1.4) .. (10.93,3.29)   ;
%Curve Lines [id:da3666427806786574] 
\draw    (211.8,66.8) .. controls (216.75,81.65) and (229.44,97.48) .. (185.26,112.54) ;
\draw [shift={(183.9,113)}, rotate = 341.64] [color={rgb, 255:red, 0; green, 0; blue, 0 }  ][line width=0.75]    (10.93,-3.29) .. controls (6.95,-1.4) and (3.31,-0.3) .. (0,0) .. controls (3.31,0.3) and (6.95,1.4) .. (10.93,3.29)   ;
%Curve Lines [id:da2239431794815021] 
\draw    (404.8,69.8) .. controls (409.75,84.65) and (422.44,100.48) .. (378.26,115.54) ;
\draw [shift={(376.9,116)}, rotate = 341.64] [color={rgb, 255:red, 0; green, 0; blue, 0 }  ][line width=0.75]    (10.93,-3.29) .. controls (6.95,-1.4) and (3.31,-0.3) .. (0,0) .. controls (3.31,0.3) and (6.95,1.4) .. (10.93,3.29)   ;
% Text Node
\draw (78,66.4) node [anchor=north west][inner sep=0.75pt]  [font=\Large]  {$p_{i}{}_{n}$};
% Text Node
\draw (550,60.4) node [anchor=north west][inner sep=0.75pt]  [font=\Large]  {$p_{o}{}_{u}{}_{t}$};
% Text Node
\draw (141,159) node [anchor=north west][inner sep=0.75pt]   [align=left] {{\large Flow front}};
% Text Node
\draw (185,43) node [anchor=north west][inner sep=0.75pt]   [align=left] {fluid 2};
% Text Node
\draw (379,46) node [anchor=north west][inner sep=0.75pt]   [align=left] {fluid 1};
\end{tikzpicture}
    \caption{One-dimensional domain (filling problem)}
    \label{fig:1d_domain}
\end{figure}

We consider a one-dimensional problem shown in figure~\ref{fig:1d_domain}. At $t=0$, the domain is filled with one fluid (fluid 1). Another fluid (fluid 2) is being injected from the left end at constant pressure $p_{in}$, while the pressure at the other end is fixed to $p_{out}$.

Equation (\ref{eq:darcy_full}) will reduce to
\begin{equation}
    v = -\dfrac{k}{\mu}p_x
    \label{eq:p_eq}
\end{equation}

where $v$ is the volume average Darcy's velocity, $k$ the permeability, $\mu$ the viscosity and $p_x$ the pressure gradient with respect to $x$. Equation~(\ref{eq:momentum_full}) reduces to
\begin{equation}
   v_x = 0
    \label{eq:v_eq}
\end{equation}

Boundary conditions are given for the pressure inlet and outlet at all times.
\begin{equation}
    p(0,t) = p_{in},\ \ p(l,t) = p_{out}
\end{equation}

where $l$ is the domain length.

To differentiate between the two phases, a fraction function $c$ is defined. $c$ takes a value of 1 in the domain that is filled with fluid~2 and 0 when fluid~1 is present. The viscosity $\mu$ is redefined according to (\ref{eq:mu}). Equation~(\ref{eq:c_full}) reduces to

\begin{equation}
       c_t + v\ c_x = 0
       \label{eq:c_eq}
\end{equation}

where $c_t$ and $c_x$ are time and position gradients of $c$, respectively.

To solve this equation, initial and boundary conditions need to be enforced on $c$.

\begin{equation}
    c(0,t) = 1,\ \ c(x, 0) = 0
\end{equation}

To sum up, the strong form of the problem can be written as:
\begin{equation}
    \begin{gathered}
    c_t + v\ c_x = 0,\ \ \ x \in [0,l],\ \ t\in [0, T],\ \\
    v = \ -k/ \mu\ p_x,\ \ \ \ \ \ \ x \in [0,l],\ \ t\in [0, T],\ \\
    v_x = 0,\ \ \ \ \ \ x \in [0,l],\ \ t\in [0, T],\ \\
    \mu = c\mu_2 + (1-c)\mu_1\\
    p(0,t) = p_{in},\ \ p(l,t) = p_{out}\\
    c(0,t) = 1,\ \ c(x, 0) = 0
\end{gathered}
\end{equation}

where $T$ is the end of the time domain where the problem is solved. 

The training is performed using a neural network of 5 hidden layers and 20 neurons in each layer. Hyperbolic tangent activation function is used in all the hidden layers; this choice proved to provide good results \cite{solid_pinn, composite_pinn, subsurface_pinn}. However, for the output layer sigmoid function is used for the pressure and fraction function networks since their values go from 0 to 1, while, linear activation function is used for the output layer of the velocity network. Adam optimizer is used for the initial training phase with a fixed learning rate of 0.001 and 1,000 iterations. Afterwards, BFGS optimization algorithm is used. This optimization strategy is followed by \cite{high_speed_pinn, res_adapt_pinn, biot_pinn}. The reason is that second-order methods like BFGS are prone to fall into local minima, therefore, Adam is firstly used to reach the zone of the global minimum, afterwards, BFGS is used to reach the minimum easier since it iteratively builds an approximation of the Hessian matrix. The different terms of the loss function are weighted similarly, thus, $\lambda_c=\lambda_p=\lambda_1=\lambda_2=\lambda_3=1$. The parameters used to solve the problem are given in table~\ref{t1}.

\begin{table}[h]
\centering
\caption{Parameters used for the one-dimensional two-phase flow.}
\begin{tabular}{|c|c|}
	\hline
	Parameter & Value \\
	\hline\hline
	$l$ & $1$\\
	\hline
	$T$ & $0.5$\\
	\hline
	$k$ & $1$\\
	\hline
	$\mu_2$ & $1$\\
	\hline
	$\mu_1$ & $10^{-5}$\\
	\hline
	$p_{in}$ & $1$\\
	\hline
	$p_{out}$ & $0$\\
	\hline
\end{tabular}
\label{t1}
\end{table}

Three numerical experiments are performed and compared. The first one using fixed number and position of collocation points (2500 points organized in 50$\times$50 grid) for the whole training phase. The second experiment using the RAR technique starting with 1600 points organized as 40$\times$40 grid points during the Adam training phase. Afterwards, point enrichments are performed every 50 BFGS iterations till the stopping criteria are satisfied (at 2500 points as well). The final experiment using the provided adaptivity algorithm starting with 1600 points organized as 40$\times$40 grid points during the Adam training phase. Similar enrichments are done as in the RAR experiment. All cases took near 200 seconds to converge using laptop Intel core i7-6700HQ CPU @ 2.60 GHz 2.59 GHz with 8 Go RAM.

A fixed 1000 points randomly distributed over the domain are used as a test set to have a sense of the generalization error committed during the training phase. It should be noted that these points are only used for testing but are not used in the training phase.

\begin{figure}[H]
    \centering
    \includegraphics[width=0.9\textwidth]{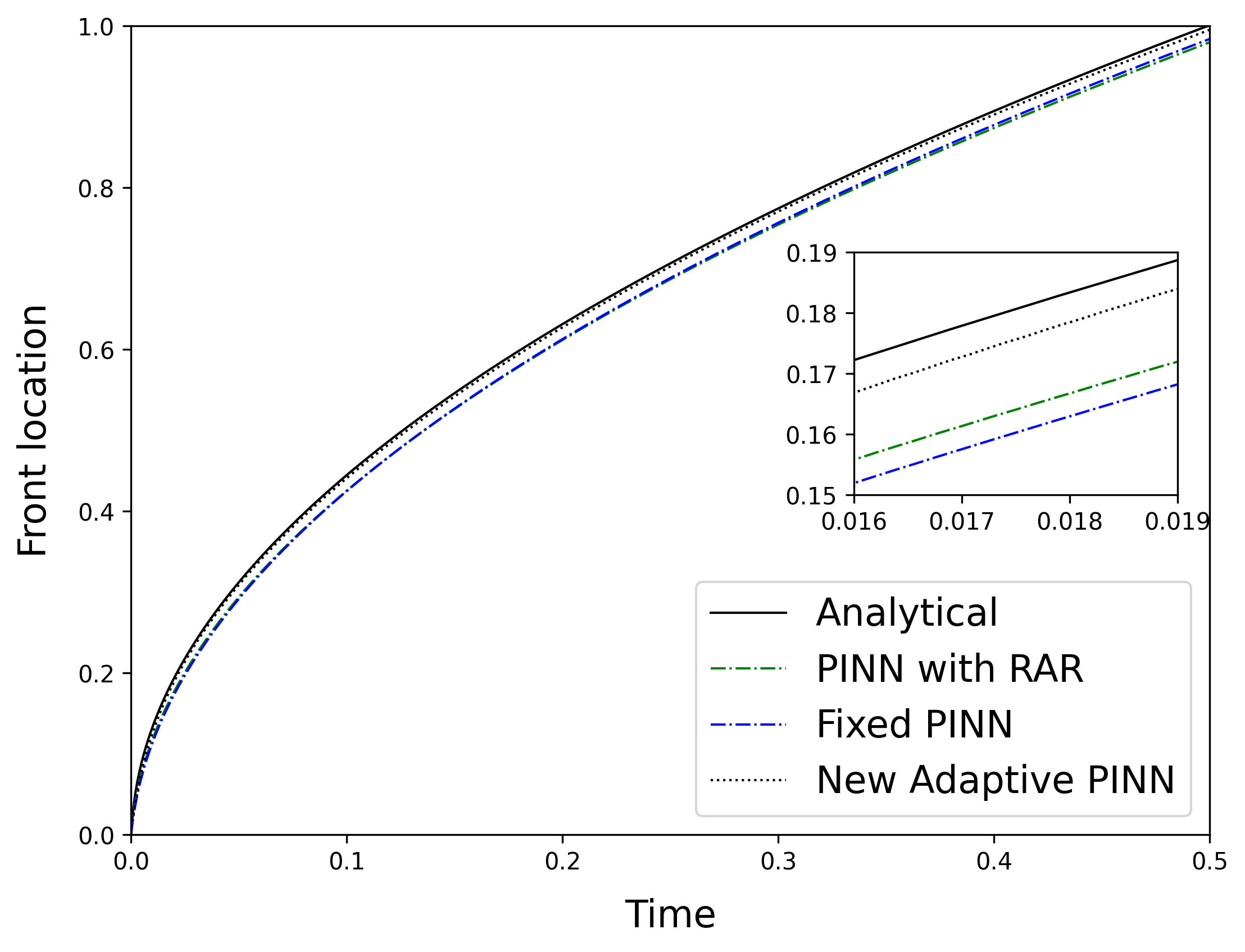}
  \caption{Front (fluid 1/fluid 2 interface) position with time for the fixed, RAR and new adaptivity cases along with analytical front position for the 1D example and a zoomed image to differentiate.}\label{fig:xf}
    
\end{figure}

Firstly, the flow front (fluid 1/fluid 2 interface) is extracted as the 0.5 level set of $c$ and plotted for both cases along with the analytical front solution in figure~\ref{fig:xf}. The analytical solution for the flow front $x_f$ is obtained as

\begin{equation}
    x_f = \dfrac{-\mu_1 l + \sqrt{\mu_1^2 l^2 + 2 (\mu_2-\mu_1)k(p_{in}-p_{out})t}}{\mu_2-\mu_1}
\end{equation}

The new adaptivity technique provided the best results among the three cases. PINN with RAR did not provide a significant improvement to the approximation; that is probably due to focusing of the enrichments in a small spatio-temporal region leading to harder optimization and higher generalization error.

\begin{figure}[H]
    \centering
    \includegraphics[width=0.9\textwidth]{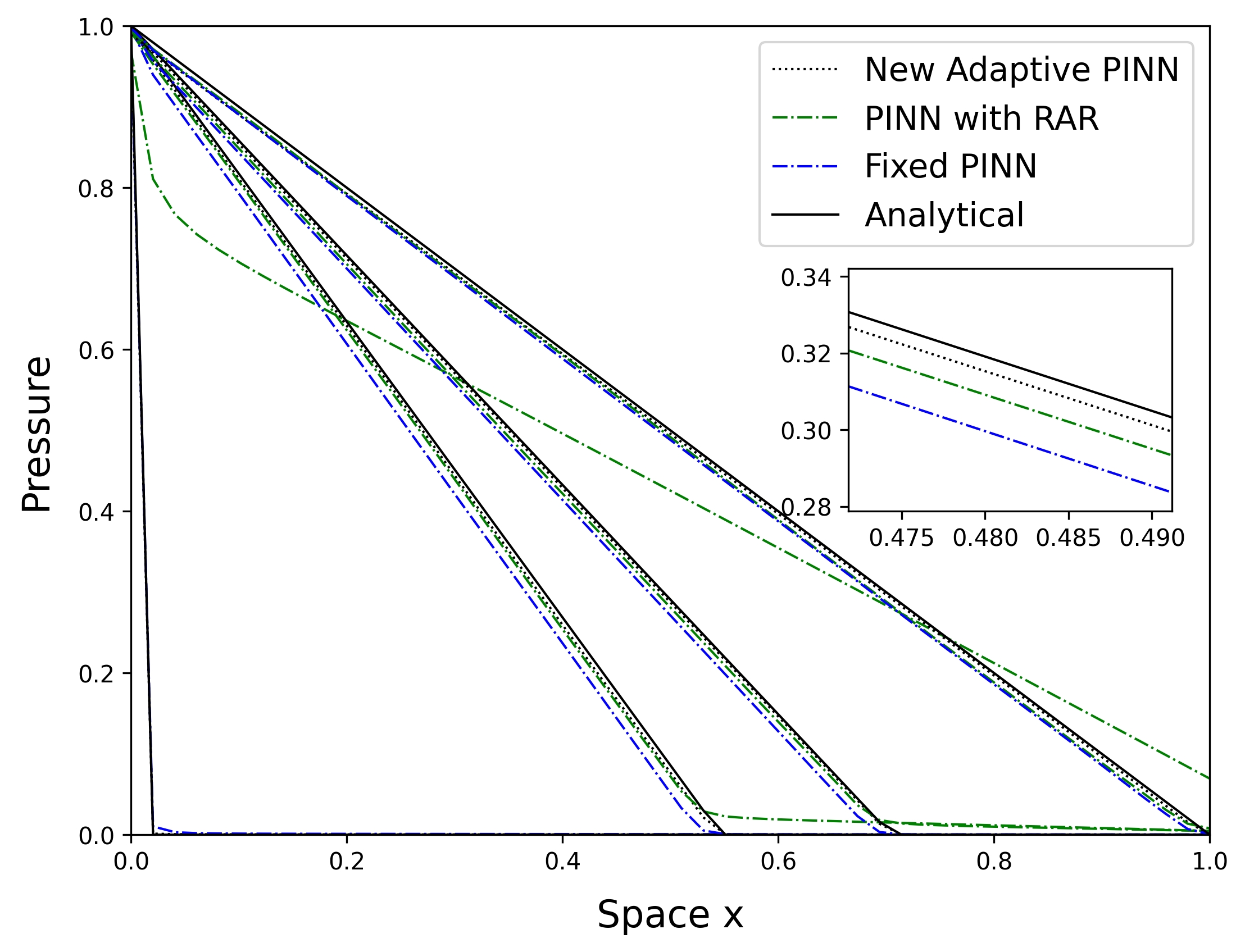}
  \caption{Pressure profiles at different times ($t=0$, $t=0.15$, $t=0.25$ and $t=0.5$) for the adaptive and fixed collocations cases for the 1D example.}\label{fig:presevol}
    
\end{figure}

The pressure profiles at different times for the three cases are shown in figure~\ref{fig:presevol} and compared to the analytical solution. The analytical solution for pressure can be written as

\begin{equation}
p(x,t)=   \left\{
  \begin{array}{ll}
      \dfrac{-\mu_2 (p_{in}-p_{out})}{(\mu_2-\mu_1)x_f(t) + \mu_1 l}\ x + p_{in} & x < x_f(t) \\[3ex]
      \dfrac{-\mu_1 (p_{in}-p_{out})}{(\mu_2-\mu_1)x_f(t) + \mu_1 l}\ x + \dfrac{\mu_1 (p_{in}-p_{out})}{(\mu_2-\mu_1)x_f(t) + \mu_1 l}\ l + p_{out} & x \geq x_f(t) \\
\end{array} 
\right.
\end{equation}

The new adaptive case provided a pressure solution closer to the analytical solution than both of the fixed PINN and PINN with RAR cases. PINN with RAR shows a prediction far from the analytical solution near $t=0$.

\begin{figure}[H]
\vspace*{-1.1in}
\hspace*{-1.25in}
    \centering
    \includegraphics[width=1.6\textwidth]{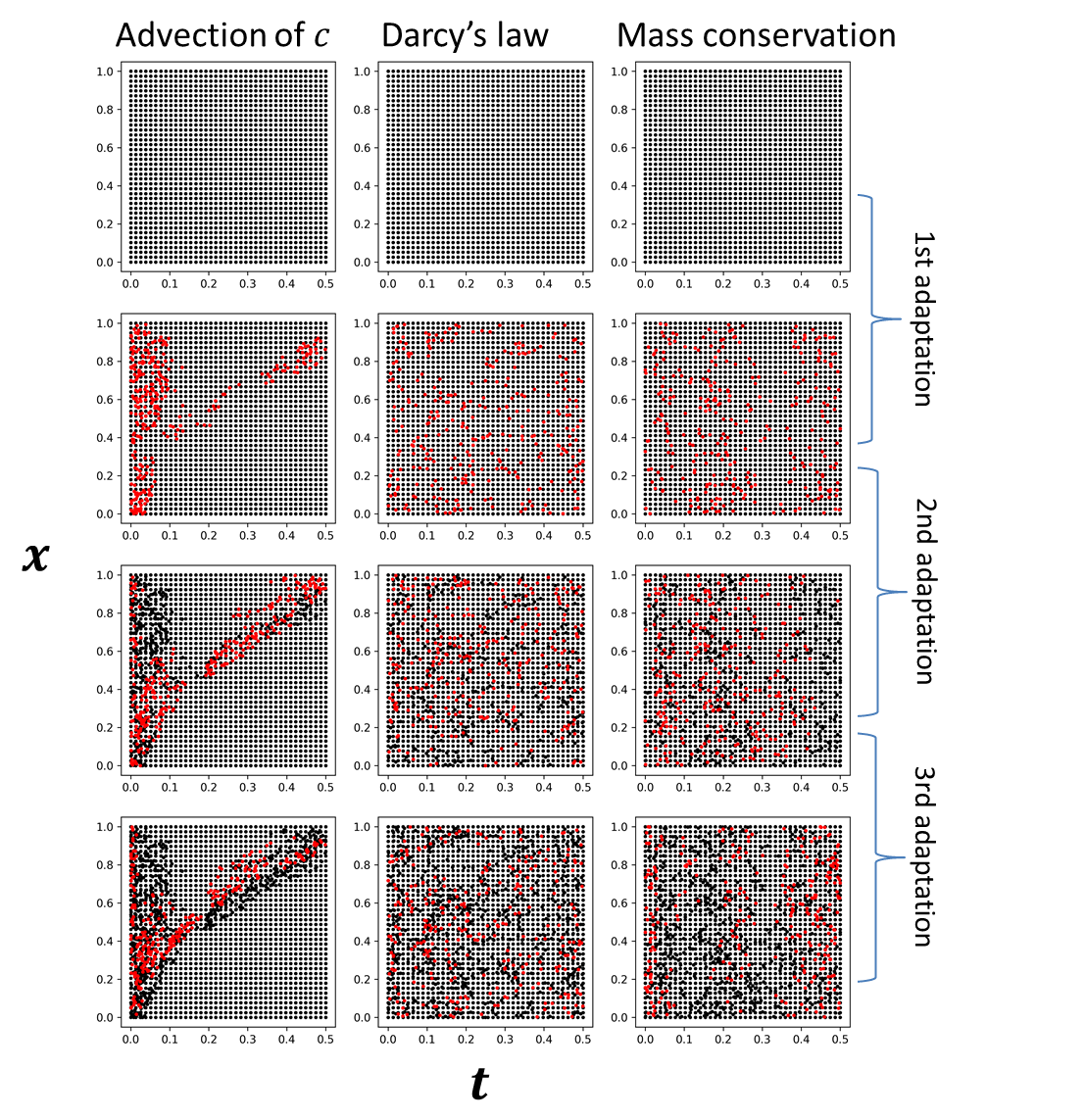}
  \caption{Evolution of collocation points with new adaptivity (from left to right: collocations for $f_1$, collocations for $f_2$ and collocations for $f_3$) for the 1D example. New collocation points after each adaptation step are shown in red.}
    \label{fig:points_collocotions}
\end{figure}

The distribution of the collocation points for the different PDEs is shown in the case of new adaptive PINN in figure~\ref{fig:points_collocotions}. The figure shows the evolution of these distributions at different stages of using the adaptivity algorithm.

From figure~\ref{fig:points_collocotions}, it can be seen for the collocation points of $f_1$ that the points seem to be dense near the location of the front, where the residual is higher. This helped in capturing the interface location accurately. While for the other collocation points ($f_2$ and $f_3$), they were distributed almost randomly in the domain. That is because the residual field is spread all over the domain since there are no sharp solution features to capture.

\begin{figure}[H]
    \centering
    \includegraphics[width=1.1\textwidth]{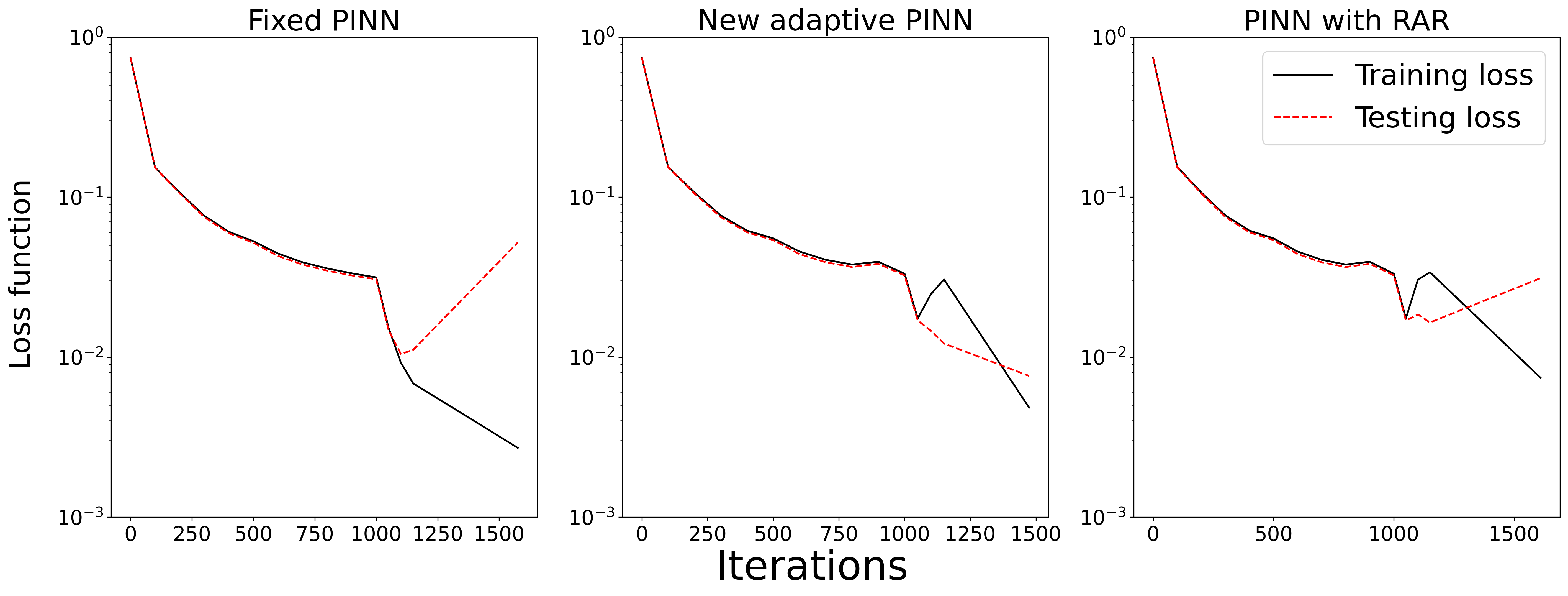}
  \caption{Loss vs. Iteration graphs for fixed PINN (left plot), new adaptive PINN (middle plot), and PINN with RAR (right plot) for the 1D example.}\label{fig:cost}
    
\end{figure}

The loss function is compared in the three cases by plotting the loss vs. iteration graph for both cases. The loss using the training set is compared to that using the testing set for all cases as shown in figure~\ref{fig:cost}. It should be noted that there is a deviation between the training and testing loss in the fixed collocations case, meaning that the generalization error is high. This deviation is marginally decreased in the PINN with RAR case. For the new adaptive PINN case, the deviation significantly decreased meaning that less generalization error is committed using the newly-developed adaptive technique. From a deep learning perspective, using adaptive collocation points can be seen as a form of regularization of the neural network solution; adaptivity prevents overfitting, thus making the solution more accurate for unseen points (points not used in the training process).

\begin{figure}[H]
    \centering
    \includegraphics[width=1\textwidth]{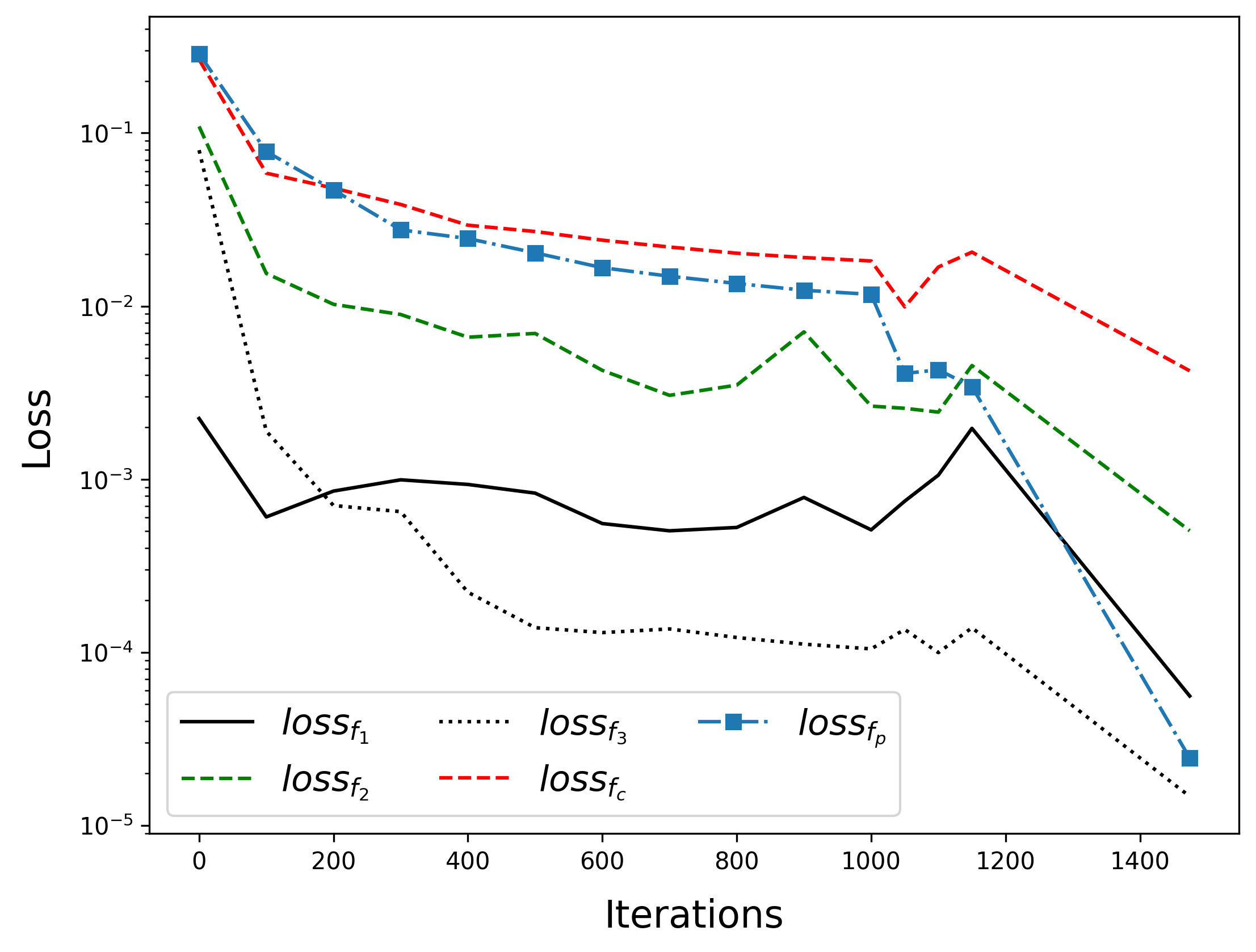}
  \caption{Different loss terms versus iterations for the new adaptivity case for the 1D example.}\label{fig:loss_terms}
\end{figure}

The different terms in the loss function are plotted for the new adaptive case against the number of iterations in figure~\ref{fig:loss_terms}. This is done to assure the convergence of all the terms.

All the loss terms have more or less a similar trend in convergence. We can note that the loss corresponding to the initial and boundary conditions of the fraction function $c$ has the highest values (harder convergence) that is due to the discontinuity in the initial/boundary condition values at $t=0$.

\subsection{Two dimensional central injection}\label{subsec:2DFilling}

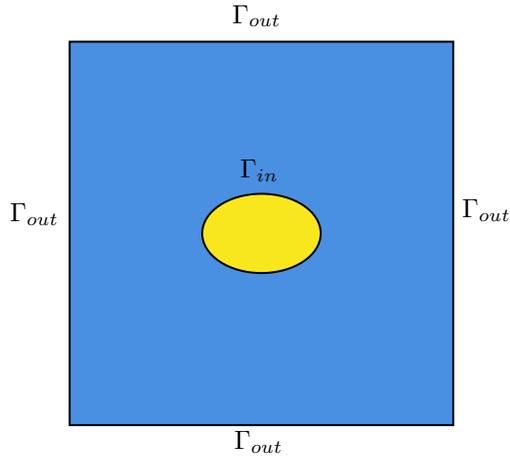
\begin{figure}[H]
\centering
\tikzset{every picture/.style={line width=0.75pt}} %set default line width to 0.75pt        

\begin{tikzpicture}[x=0.75pt,y=0.75pt,yscale=-1,xscale=1]

%Shape: Square [id:dp689705754189113] 
\draw  [fill={rgb, 255:red, 74; green, 144; blue, 226 }  ,fill opacity=1 ] (164.47,52.13) -- (358,52.13) -- (358,245.67) -- (164.47,245.67) -- cycle ;
%Shape: Ellipse [id:dp46068327458118774] 
\draw  [fill={rgb, 255:red, 248; green, 231; blue, 28 }  ,fill opacity=1 ] (231.35,148.9) .. controls (231.35,137.85) and (244.73,128.9) .. (261.23,128.9) .. controls (277.74,128.9) and (291.12,137.85) .. (291.12,148.9) .. controls (291.12,159.95) and (277.74,168.9) .. (261.23,168.9) .. controls (244.73,168.9) and (231.35,159.95) .. (231.35,148.9) -- cycle ;

% Text Node
\draw (245,32) node [anchor=north west][inner sep=0.75pt]   [align=left] {$\displaystyle \Gamma _{out}$};
% Text Node
\draw (246,247) node [anchor=north west][inner sep=0.75pt]   [align=left] {$\displaystyle \Gamma _{out}$};
% Text Node
\draw (133,132) node [anchor=north west][inner sep=0.75pt]   [align=left] {$\displaystyle \Gamma _{out}$};
% Text Node
\draw (361,130) node [anchor=north west][inner sep=0.75pt]   [align=left] {$\displaystyle \Gamma _{out}$};
% Text Node
\draw (249,110) node [anchor=north west][inner sep=0.75pt]   [align=left] {$\displaystyle \Gamma _{in}$};

\end{tikzpicture}
\caption{Domain of the 2D central injection problem. $\Gamma_{out}$ is the outlet boundary where the pressure is set to 0, while $\Gamma_{in}$ is the boundary of the inlet where pressure is set to 1.}
    \label{fig:2D_domain}
\end{figure}

The next example is a two-dimensional problem that is encountered in structural composites manufacturing processes such as resin transfer molding. The domain is a square of a unity area with an elliptic injection port placed at the center with constant pressure ($p=1$). The four outer sides are outlets where the pressure is set to zero. The problem domain is plotted in figure~\ref{fig:2D_domain}. The analytical solution of this problem exists in \cite{analytical_2d}.

 The material properties are the same as in the first example except for the permeability of the domain which is uniform orthotropic and can written as:

\begin{equation}
   \mathbf{K} = \begin{bmatrix}
1.5 & 0 \\
0 & 1 
\end{bmatrix}
\end{equation}

All training parameters, that are used in the first example, are the same in this example except for the network architecture. In this case, 4 neural networks are used to approximate the fraction function, pressure and two velocity components; each network is composed of 5 hidden layers and 20 neurons.

\begin{figure}[H]
    \centering
    \includegraphics[width=1.1\textwidth]{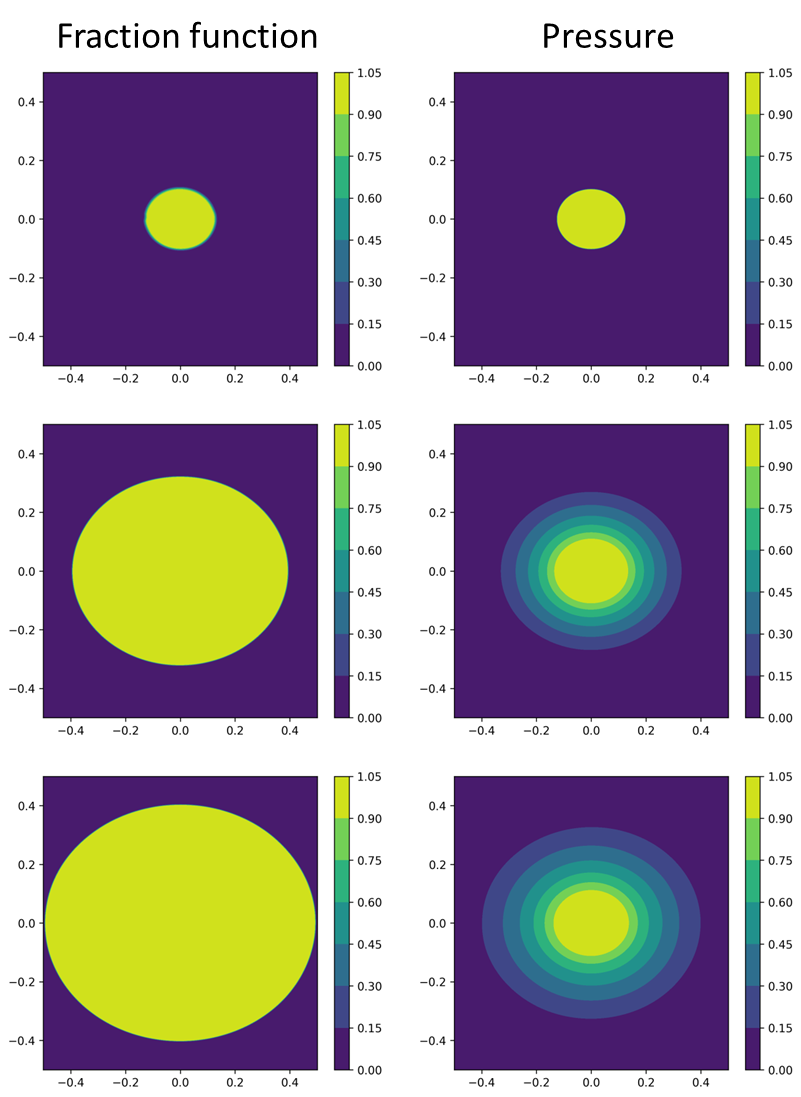}
  \caption{Left: Front function field. Right: pressure field. The plots are shown at times 0, 0.039 and 0.078 from top to bottom.}
    \label{fig:flow_visualize}
\end{figure}

The fraction function and pressure fields are plotted at different times for the adaptive case to visualize the evolution of the flow front and pressure with time in figure~\ref{fig:flow_visualize}.

\begin{figure}[H]
  \centering
  \includegraphics[width=1\textwidth]{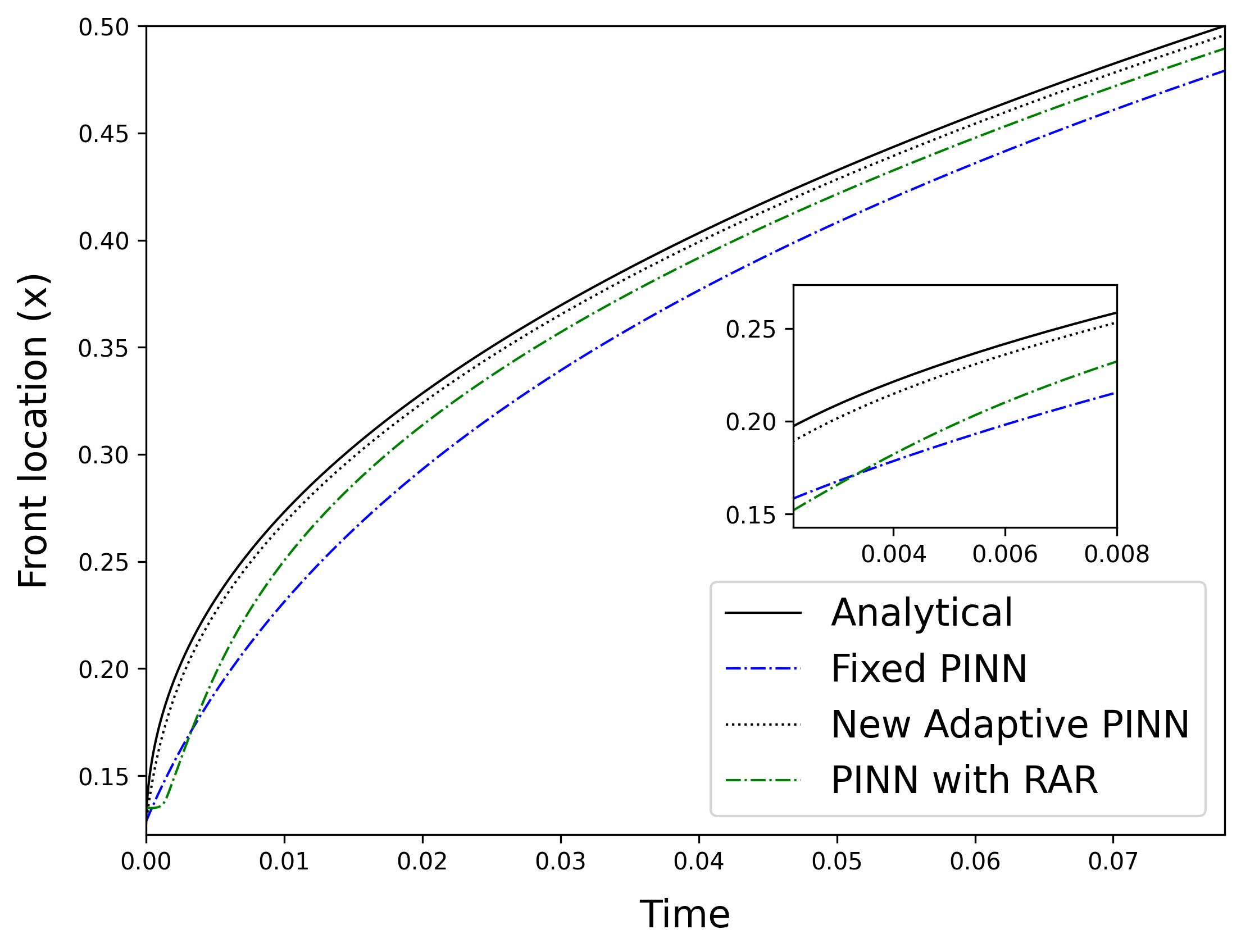}
  \caption{Front (fluid 1/fluid 2 interface) position in the x-direction with time for the fixed and adaptive cases along with analytical front position for the 2D case.}
  \label{fig:rf_x}
\end{figure}

\begin{figure}[H]
  \centering
  \includegraphics[width=1\textwidth]{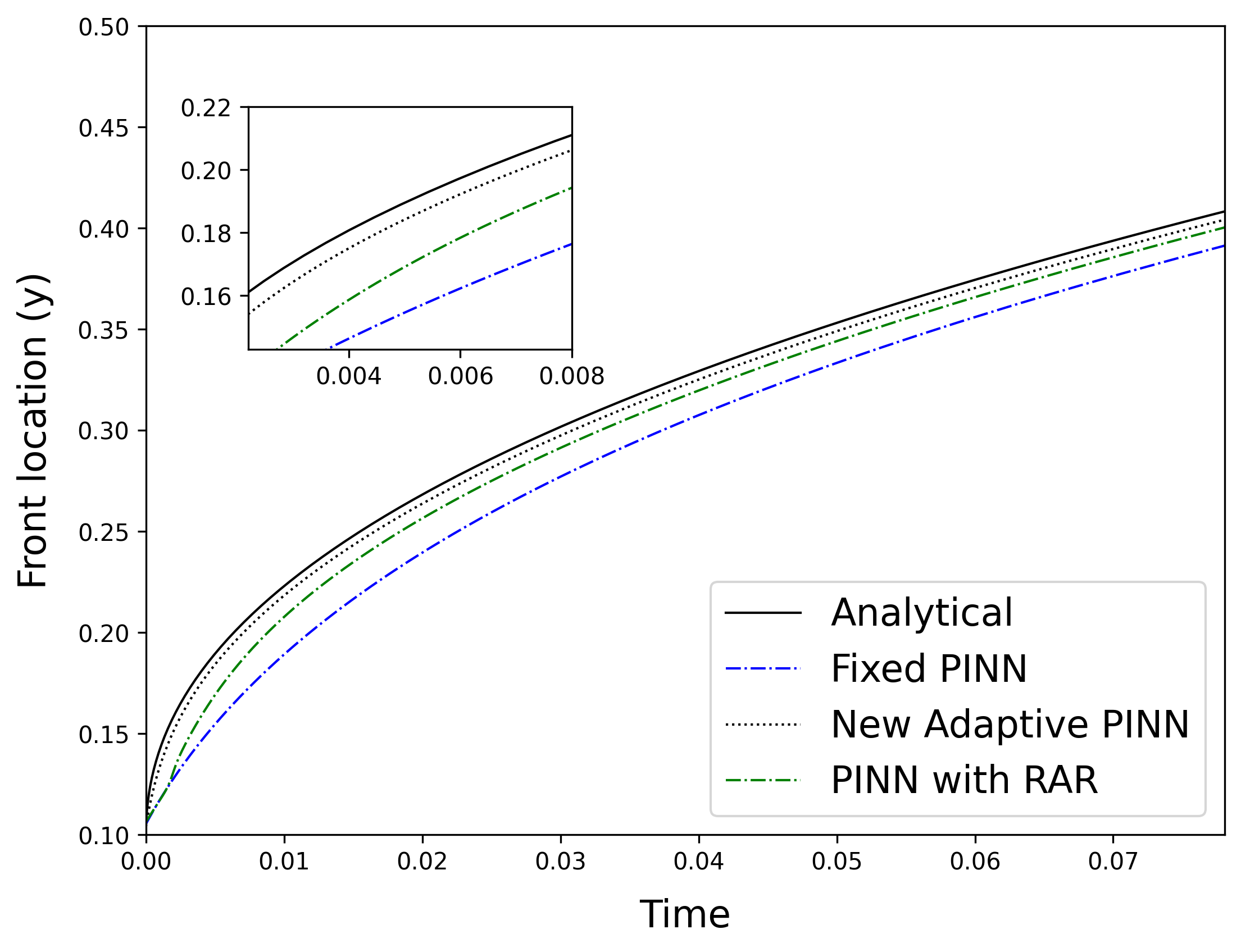}
  \caption{Front (fluid 1/fluid 2 interface) position in the y-direction with time for the fixed and adaptive cases along with analytical front position for the 2D case.}
  \label{fig:rf_y}
\end{figure}

The flow front position in x and y directions is plotted in figures~\ref{fig:rf_x} and \ref{fig:rf_y} for the 3 cases along with the analytical solution to compare.

It can be seen that the new adaptive algorithm provided better prediction of the flow front position.

\begin{figure}[H]
    \centering
    \includegraphics[width=1.1\textwidth]{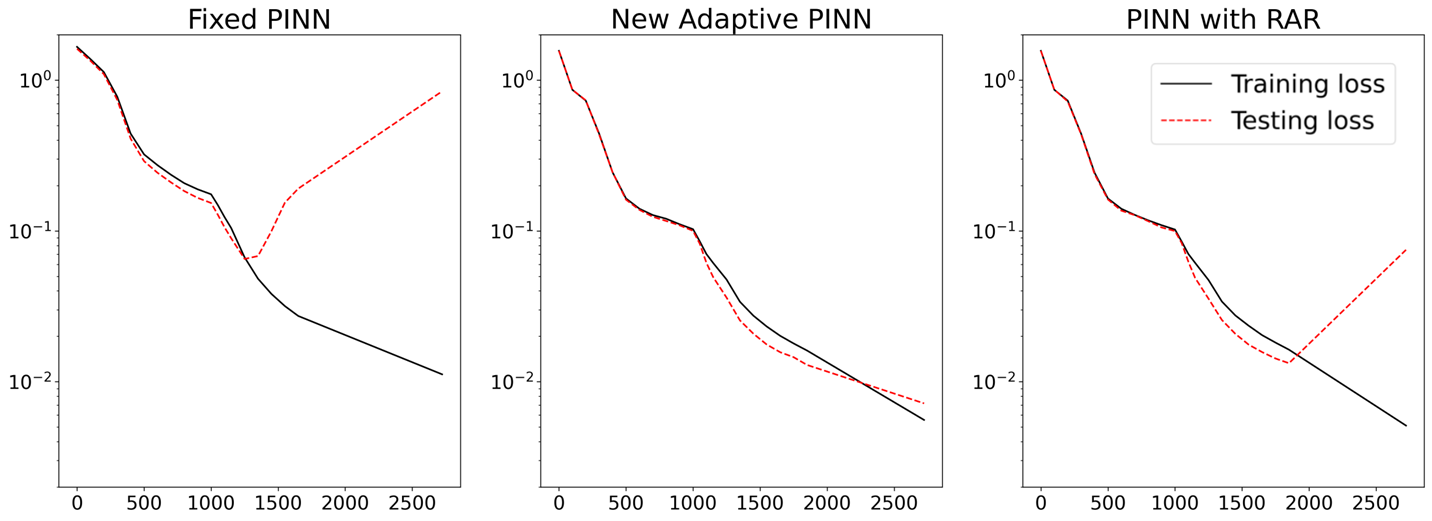}
  \caption{Loss vs. Iteration graphs for fixed PINN (left plot), new adaptive PINN (middle plot), and PINN with RAR (right plot).}
    \label{fig:2d_loss}
\end{figure}

The loss function is plotted for the 3 cases in figure~\ref{fig:2d_loss}. 10,000 points in the space-time domain are chosen randomly to asses the testing loss, while they are not used in the training.

Using the new adaptivity algorithm, the discrepancy between the training and testing loss is greatly reduced. That means that it offers a mean to reduce the generalization error thus provides better accuracy to the approximated fields.

\section{Conclusion}
\label{section:conclusion}

In this study, Physics-Informed Neural Networks (PINN) are applied to solve two-phase flow problems in porous media. A novel residual-based adaptive algorithm is developed. The key is to utilize the PDE residual to build a probability density function, from which new collocation points are drawn and added to the training set. The technique is applied to the different PDEs in the coupled system independently, thus, different collocations points are used for the different PDEs. Moreover, the technique is applied to enrich the points used to capture the initial and boundary conditions. 

The adaptivity algorithm is used to solve 1D and 2D two-phase flow in porous media . Using the new technique provided better results than using the classical PINN with fixed collocation points and also showed an improvement over the RAR technique. The adaptive technique can be seen as a form of regularization of the neural networks, thus reducing the generalization error.

\section*{Acknowledgements}

This study was funded under the PERFORM Thesis program of IRT Jules Verne, Bouguenais, France.

\bibliography{biblio.bib}

\end{document}